\documentstyle[12pt,amscd]{amsart}
\newtheorem{Theorem}{Theorem}[section] \newtheorem{Lemma}[Theorem]{Lemma}
\newtheorem{Corollary}[Theorem]{Corollary}
\newtheorem{Proposition}[Theorem]{Proposition}

\def\height{\operatorname{ht}} 
\def\reg{\operatorname{reg}} \def\ini{\operatorname{in}}
\def\gin{\operatorname{Gin}} \def\Min{\operatorname{Min}}
 \def\To{\longrightarrow}

\def\sat{{sat}}

\def\To{\longrightarrow}
\def\sk{\smallskip}
\def\mm{{\frak m}} \def\nn{{\frak n}}
\def\pp{{\frak p}} \def\qq{{\frak q}}

\begin{document}
\title{Asymptotic linear bounds for the Castelnuovo-Mumford regularity}
\author{J\"urgen Herzog, L\^e Tu\^an Hoa and Ng\^o Vi\^et Trung}
\address{Fachbereich Mathematik,  Universit\"at-GHS Essen, 45117 Essen, Germany}
\email{juergen.herzog@@uni-essen.de}
\address{Institute of Mathematics,  Box 631, B\`o H\^o, 10000 Hanoi, Vietnam}
\email{lthoa@@hanimath.ac.vn}
\address{Institute of Mathematics,  Box 631, B\`o H\^o, 10000 Hanoi, Vietnam}
\email{nvtrung@@hn.vnn.vn}
\thanks{The second and third authors are partially supported by the National Basic Research 
Program of Vietnam} 
\begin{abstract} We prove asymptotic linear bounds for the Castelnuovo-Mumford regularity of 
certain filtrations of homogeneous ideals whose 
Rees algebras need not to be Noetherian. \end{abstract} 
\maketitle

\section*{Introduction}

In a pioneering paper, Bertram, Ein and Lazarsfeld [BEL] proved that if $X \subset {\Bbb P}^r$ 
is a smooth complex variety of codimension $s$ which is cut out  scheme-theoretically by 
hypersurfaces of degree $d_1 \ge \ldots \ge d_m$, then $H^i({\Bbb P}^r,{\cal I}_X^n(a)) = 0$ 
for $i \ge 1$ and $a \le d_1n + d_2 + \cdots + d_s - r$. Their result has initiated the study 
on the Castelnuovo-Mumford regularity of the powers of a homogeneous ideal. The first result 
in this direction is due to Chandler [Ch], Geramita, Gimigliano and Pitteloud [GGP]. They 
proved that if $I$ is a homogeneous ideal in a polynomial ring $R$ (over a field) with $\dim 
R/I \le 1$, then $\reg(I^n) \le \reg(I)n$ for all $n \ge 1$. This result does not hold for 
higher dimension due to an example of Sturmfels [St]. However, Swanson [Sw] showed that 
$\reg(I^n)$ is always bounded by a linear function of the form $cn+e$. She could not provide 
an estimate for $c$. For monomial ideals, there were some attempts to estimate $c$ in terms of 
better understood invariants of $I$ (see [SS], [HT]). This problem seemed to be hard. So it 
came as a surprise that we can always choose $c \le d(I)$, where $d(I)$ denotes the maximum 
degree of the generators of $I$, and that $\reg(I^n)$ is in fact a linear function for all 
large $n$. \sk

\begin{Theorem} \label{ordinary} {\rm  [CHT, Theorem 3.1], [Ko, Theorem 5]} Let $I$ be an 
arbitrary homogeneous ideal in a polynomial ring. There exist constants $c \le d(I)$ and $e 
\ge 0$ such that for all large $n$,
$$\reg(I^n) = cn + e.$$ \end{Theorem} \sk

The reason for this phenomenon comes from the fact that the Rees
algebra $R[It] = \oplus_{n \ge 0}I^nt^n$ is a Noetherian bigraded
algebra. In fact, it can be shown that if $E = \oplus_{n \ge
0}E_n$ is a finitely graded bigraded module over the Rees algebra
$R[It]$, then $\reg(E_n)$ is a linear function of $n$ for all
large $n$ (see [CHT, Lemma 3.3 and Theorem 3.4]). \sk

Recently, Cutkosky, Ein and Lazarsfeld proved the following
related, more geometric result. \sk

\begin{Theorem} \label{saturation} {\rm  [CEL, Theorem 3.2] } 
Let $I$ be an arbitrary homogeneous ideal in a polynomial ring. Let
$(I^n)^\sat$ denote the saturation of $I^n$. Then $\displaystyle
\lim_n\frac{\reg((I^n)^\sat)}{n}$ exists.
\end{Theorem}

This result is somewhat surprising since the saturated Rees algebra
$\oplus_{n \ge 0}(I^n)^\sat t^n$ needs not to be Noetherian.
Moreover, the limit can be an irrational number [Cu].\sk

Inspired of the above results we raise the following questions:\sk

\noindent{\bf Question 1.} Let $(R,\mm)$ be a local ring and $I$
an ideal in $R$. Let $G(R/I^n)$ denote the associated graded ring
of $R/I^n$ with respect to $\mm$. Then $G(R/I^n)$ can be represented as a quotient ring of the 
associated graded ring $G(R)$ by a homogeneous ideal $I^*_n$, where $\{I^*_n\}$ is a 
filtration of ideals. Does $\displaystyle \lim_{n\to \infty} \frac{\reg(G(R)/I_n^*)}{n}$ 
exist? \sk

\noindent{\bf Question 2.} Let $R$ be a polynomial ring and $I$ a
homogeneous ideal in $R$. Let $I^{(n)}$ denote the $n$th symbolic
power of $I$. Does $\displaystyle \lim_{n\to
\infty}\frac{\reg(I^{(n)})}{n}$ exist? \sk

\noindent{\bf Question 3.} Let $R$ be a polynomial ring and $I$ a
homogeneous ideal in $R$. Let $\ini(I)$ denote the intial ideal of
$I$ with respect to an arbitrary term order. Does $\displaystyle
\lim_{n\to \infty}\frac{\reg(\ini(I^n))}{n}$ exist? \sk

In all these cases, the underlying bigraded ``Rees algebra" need not to be Noetherian, hence 
the method of [CHT] and [Ko] could not be
applied directly. In fact, we do not know any general approach and any definite answer (yes or 
no) to these problems. \sk

One may weaken the above questions by asking whether the Castelnuovo-Mumford regularities 
$\reg(G(R/I_n))$, $\reg(I^{(n)})$ and $\reg(\ini(I^n))$ are bounded by linear functions of 
$n$. In this paper we will prove the following results which suggest that the weakened 
questions may have a positive answer: \sk

\noindent 1. $\reg(G(R/I^n))$ is bounded by a linear function if $\dim R/I \le 1$.\sk

\noindent 2. $\reg(I^{(n)})$ is bounded by a linear function in the follwong cases:\par
(i) $\dim R/I \le 2$ (already discovered by Chandler [Ch]),\par
(ii) The singular locus of $R/I$ has dimension $\le 1$,\par
(iii) $I$ is an arbitrary monomial ideal. \sk

\noindent 3. $\reg(\ini(I^n)) \le \reg(\ini(I))n$ and $\displaystyle \lim_{n\to \infty} 
\frac{\reg(\ini(I^n))}{n}$ exits if $\dim R/I \le 1$. \sk

Moreover, the $a$-invariant of the graded rings $G(R/I^n)$, $R/I^{(n)}$ and $R/\ini(I^n)$ is 
bounded by a linear function of $n$ for arbitrary dimension. Note that the Castelnuovo-Mumford 
regularity is always an upper bound for the $a$-invariant. \sk

The paper is divided in three sections. \sk

In Section 1 we first discuss the relationship between the Castelnuovo-Mumford regularity and 
the reduction number which is also a measure for the complexity of a graded ring. Then we 
apply this relationship to study Question 1. Let $r(R/I^n)$ denote the reduction number of the 
maximal ideal of $R/I^n$. We show that $c = \displaystyle \lim_{n\to \infty} 
\frac{r(R/I^n)}{n}$ exits and that $r(R/I^n) \le cn-1$ for all $n$. The answer to Problem 1 in 
the case $\dim R/I \le 1$ follows from this fact and Swanson's linear bound for the nilpotency 
index of primary components of $I^n$ (see [Sw]). Moreover, we construct examples with 
$\reg(G(R/I^n)) = r(R/I^n)$ for all $n \ge 1$ such that $\displaystyle \lim_{n\to \infty} 
\frac{r(R/I^n)}{n}$ is any given rational number $\ge 1$. \sk

In Section 2 we study the Castelnuovo-Mumford regularity of the generalized symbolic powers 
which are defined by taking out those primary components of $I^n$ which contain a given ideal 
$J$. If $\dim R/J \le 2$, we can bound $\reg(I^{(n)})$ by means of $\reg(I^n)$. The first two 
cases of the above-mentioned results to Question 2 are simple consequences of this fact. The 
case of monomial ideals follows from the fact that the regularity of a monomial ideal is 
bounded by the degree of the least common multiple of the generators (see [BH2] and [HT]). \sk

In Section 3 we study the initial ideal $\ini_{\lambda}(I)$ with respect to an integral weight 
function $\lambda$ instead of the usual initial ideal $\ini(I)$ with respect to a term order. 
We show that the Castelnuovo-Mumford regularity and the reduction number of the quotient ring 
$R/\ini_{\lambda}(I)$ do not change when passing to the generic member of certain flat family 
over $k[t]$ of quotient rings of $S$. This relationship allows us to prove the above-mentioned 
result to Question 3.
As a byproduct of our approach, we can apply a recent result of Ein, Lazarsfeld and Smith 
[ELS] to show that for all $n \ge 0$, $\ini(I^{sn}) \subseteq \ini(I)^n$, where $s$ is the 
codimension of $I$. \sk

\noindent{\it Acknowledgement.} This paper was written
during a visit of L\^e Tu\^an Hoa and Ng\^o Vi\^et Trung to the
University of Essen im Summer 2000. These authors would like to
express their sincere thanks to the Research Group "Arithmetic and
Geometry" for generous supports and hospitality.\sk

\section{Regularity and reduction number} \sk

Let us first recall some facts on the Castelnuovo-Mumford
regularity of a graded module and its relationship to the
reduction number. We will do it in a general setting (since it
will be needed later in this paper). \sk

Let $S = S_0[x_1,\ldots,x_m]$ be a standard graded algebra over a
local ring $(S_0,\nn)$. For convenience we assume that the residue
field of $S_0$ is infinite. Let $E$ be a finitely graded
$S$-module. For every integer $i \ge 0$ we set
$$a_i(E) := \max\{r|\ H_{S_+}^i(E)_r \neq 0\}$$
with $a_i(E) = -\infty$ if $H_{S_+}^i(E) = 0$, where
$H_{S_+}^i(E)$ denotes the $i$th local cohomology module of $E$
with respect to the ideal $S_+ = \oplus_{n>0}S_n$.\sk

\noindent{\bf Definition}. The {\it Castelnuovo-Mumford
regularity} of $E$ is defined by
$$\reg(E) = \max\{a_i(E)+i|\ i \ge 0\}.$$

Note that $H_{S_+}^i(E) = 0$ for $i > \dim E/\nn E$ (see e.g.
[T1], [T2]). We will call $E/\nn E$ the fiber module of $E$.\sk

\noindent{\it Remark.} Let $J$ be a homogeneous ideal of $S$.
Since $H_{S_+}^i(S) = 0$ for $i \neq m$ and $H_{S_+}^i(J) = 0$
for $i > m$, from the exact sequence $0 \To J \To S \To S/J$ we
can deduce that $\reg(J) = \reg(S/J)+1.$ \sk

\noindent{\bf Definition}. A {\it reduction} of $E$ is a graded
ideal $J$ generated by forms of degree 1 such that $(JE)_n = E_n$
for all large $n$. The least integer $n_0$ for which $(JE)_n =
E_n$ for all $n > n_0$ is called the {\it reduction number} of $E$
with respect to $J$. It will be denoted by $r_J(E)$. A reduction
$J$ of $E$ is {\it minimal} if $J$ does not contain any other
reduction of $E$. We set
$$r(E) := \min\{r_J(E)|\ \text{$J$ is a minimal reduction of
$E$}\}.$$

 These notions were originally introduced for graded
quotient rings of $S$ (see e.g. [T1], [V]). \sk

\noindent{\it Remark.} If $S$ is a polynomial ring over a field
$k$, then $J = (z_1,\ldots,z_d)$ is a minimal reduction of a
graded quotient ring $S/I$ if and only if $A = k[z_1,\ldots,z_d]$
is a Noether normalization of $S/I$. Moreover, $r_J(S/I)$ is the
maximum degree of the generators of $S/I$ as a graded $A$-module
(see e.g. [V]). \sk

The following relationship between regularity and reduction number
has been proved for graded quotient rings of $S$. But it can be
easily extended to the case of graded modules. We leave the
reader to check the proof. \sk

\begin{Proposition} \label{bound-reduction} {\rm [T1, Proposition 3.2]}
Let $d = \dim E/\nn E$. Then
$$a_d(E)+d  \le r(E) \le \reg(E).$$
\end{Proposition} \sk

If $S$ is a polynomial ring over a field, we set $a(E) = a_d(E)$
and call it the $a$-invariant of $E$. This invariant plays an
important role in local duality since $-a(E)$ is the initial
degree of the canonical module of $E$ (see e.g. [GW], [BH1]). \sk

The reduction number of a graded module is a generalization of 
the reduction number in the local case. Let $(R,\mm)$ be a local ring. 
For convenience we assume that the residue field $k = R/\mm$ is infinite. \sk

\noindent{\bf Definition.} We call an $\mm$-primary
ideal $\qq$ a reduction of $\mm$ if $\mm^{r+1}  = \qq\mm^r$. The
least integer $r \ge 0$ with this property is called the reduction
number of $R$ with respect to $\qq$. It will be denoted by
$r_\qq(R)$. A reduction of $\mm$ is called minimal if it does not
contain any other reduction of $\mm$. We will set
$$r(R) = \min\{r_\qq(R)|\ \text{$\qq$ is a minimal reduction of
$\mm$}\}.$$ 

Let $G(R)= \oplus_{n\ge 0}\mm^n/\mm^{n+1}$ be the associated
graded ring of $R$ with respect to $\mm$. Then $G(R)$ is a
standard graded $k$-algebra. There is the following relationship between the reduction numbers 
of $R$ and $G(R)$. \sk

\begin{Lemma} \label{comparison-reduction} $r(R) = r(G(R))$. \end{Lemma}

\begin{pf} For any element $z \in \mm\setminus \mm^2$ we
denote by $z^*$ the initial form of $z$ in $G(R)$. Let
$z_1,\ldots,z_d$ be arbitrary elements in $\mm \setminus \mm^2$, 
$d = \dim R$. It is known that $(\underline z) = (z_1,\ldots,z_d)$ 
is a minimal reduction of $\mm$ if and only if $(\underline z^*) 
= (z_1^*,\ldots,z_d^*)$ is a minimal reduction of $G(R)$ and that 
$r_{(\underline z)}(R) = r_{(\underline z^*)}(G(R))$ [NR]. 
Hence the conclusion is immediate. \end{pf} \sk

The regularity of the associated graded ring $G(R)$ can be 
characterized as follows.\sk

\begin{Theorem} \label{Reg2}  {\rm [T2, Theorem 1.1] }
Let $(R,\mm)$ be a local ring with $d = \dim R$. Let $\qq$ be a
minimal reduction of $\mm$. Then the following conditions are
equivalent for a fixed integer $r \ge 0$, \par  {\rm (i) }
$\reg(G(R)) = r$,\par {\rm (ii) } There is a minimal basis
$z_1,\ldots,z_d$ for $\qq$ such that
$$[(z_1,\ldots,z_{i-1}) : z_i] \cap \mm^{r+1} = (z_1,\ldots,z_{i-1})\mm^r,\
i=1,\ldots,d,$$ and $r$ is the least integer $\ge r_\qq(R)$ with
this property.
\end{Theorem} \sk

\noindent{\it Remark.} If $G(R)$ is a Cohen-Macaulay ring, then 
the initial forms of $z_1,\ldots,z_d$ in $G(R)$ form a regular sequence.
By [VV] this implies that $z_1,\ldots,z_d$ is a regular sequence of $R$ and
that $(z_1,\ldots,z_{i-1}) \cap \mm^{n+1} = (z_1,\ldots,z_{i-1})\mm^n$
for all $n$ and $i=1,\ldots,d$. Thus, $\reg(G(R)) = r_\qq(R) = r(R)$ 
for any minimal reduction $\qq$ of $\mm$ (cf. [T2, Theorem 6.4]). \sk

The above results show that the reduction number can be used to
estimate the Castelnuovo-Mumford regularity. We shall see that 
the reduction number in the local case has good asymptotic properties 
(the same result holds for the graded case and we leave the reader 
to check it). First, we can always find a generic minimal reduction 
with the smallest reduction number. \sk

\begin{Lemma} \label{generic} Let $(R,\mm)$ be a local ring
with $d = \dim R$. For a generic choice of elements $z_1,\ldots,z_d$
in $\mm \setminus \mm^2$, the ideal $\qq = (z_1,\ldots,z_d)$ is 
a minimal reduction of $\mm$ with $r_\qq(R) = r(R)$. \end{Lemma}

\begin{pf} It has been shown in [T3, Lemma 4.2] that for a generic choice
of linear forms $z_1^*,\ldots,z_d^*$ in the associated graded ring $G(R)$,
the ideal $J = (z_1^*,\ldots,z_d^*)$ is a minimal reduction of $G(R)$
with $r_J(G(R)) = r(G(R))$. If $z_1,\ldots,z_d$ are elements in 
$\mm \setminus \mm^2$ whose initial forms in $G(R)$ are $z_1^*,\ldots,z_d^*$, 
then $\qq = (z_1,\ldots,z_d)$ is also a minimal reduction of $\mm$ with 
$r_\qq(R) = r_J(G(R))$ [NR]. Hence the conclusion follows from 
Lemma \ref{comparison-reduction}. \end{pf} \sk

\begin{Theorem} \label{asymptotic-reduction} Let $(R,\mm)$ be a
local ring. Let $\{I_n\}$ be a filtration of ideals in $R$ with
$\dim R/I_n = \dim R/I_1$ for all $n \ge 0$. Then\par
 {\rm (i) } $r(R/I_n) \le [r(R/I_1)+1]n-1$ for all $n \ge 0$,
\par {\rm (ii) } $c = \displaystyle \lim_{n \to \infty} \frac{r(R/I_n)}{n}$ exists
and $r(R/I_n) \ge cn-1$ for all $n \ge 0$. \end{Theorem}

\begin{pf} It suffices to show that for all integers $a, b \ge 0$,
$$r(R/I_{a+b}) \le r(R/I_a) + r(R/I_b) + 1.$$
Indeed, (i) is an immediate consequence of this formula. Further,
this formula implies that $c = \displaystyle \lim_{n \to \infty}
\frac{r(R/I_n)+1}{n}$ exists and that $r(R/I_n) + 1 \ge cn$ for
all $n \ge 0$ (see the remark below), hence (ii). \par 
To prove the above formula we let $d = \dim R/I_1$, $r =
r(R/I_a)$ and $s = r(R/I_b)$. Using Lemma
\ref{generic} we can find elements $z_1,\ldots,z_d \in \mm$ such
that $\underline z = \{z_1,\ldots,z_d\}$ generate minimal
reductions of the maximal ideals of the local rings $R/I_a$ and
$R/I_b$ with $r_{(\underline z)}(R/I_a) = r$ and $r_{(\underline
z)}(R/I_b) = s$. Since $\mm^{r+1} + I_a  =  (\underline z)\mm^r
+I_a$ and $\mm^{s+1} + I_b = (\underline z)\mm^s + I_b$ we get
\begin{eqnarray*} \mm^{r+1} & \subseteq &
(\underline z)\mm^r +I_a \cap \mm^{r+1},\\ \mm^{s+1} & \subseteq &
(\underline z)\mm^s + I_b \cap \mm^{s+1}.\end{eqnarray*} From this
it follows that
\begin{eqnarray*} \mm^{r+s+2} & \subseteq &
[(\underline z)\mm^r +I_a \cap \mm^{r+1}][(\underline z)\mm^s +
I_b \cap \mm^{s+1}]\\
& \subseteq & (\underline z)\mm^{r+s+1} + I_{a+b}.
\end{eqnarray*} Now, we can conclude that $\mm^{r+s+2} + I_{a+b}
= (\underline z)\mm^{r+s+1} + I_{a+b}$. Hence $\underline z$
generates a minimal reduction of the maximal ideal of the local
ring $R/I_{a+b}$ with $r_{(\underline z)}(R/I_{a+b}) \le r+s+1$.
Thus, $r(R/I_{a+b}) \le r + s +1 = r(R/I_a) + r(R/I_b) + 1$.
\end{pf}\sk

\noindent{\it Remark.} Let $\{c_n\}$ be any sequence of real
positive numbers with the property
$$c_{a+b} \le c_a + c_b$$
for all $a, b \ge 0$. Then $ c=\displaystyle \lim_{n\to\infty}
\frac{c_n}{n}$ exists and $c_n \ge cn$ for all $n \ge 0$. Indeed,
the sequence $\{\displaystyle\frac{c_{2^n}}{2^n}\}$ is monotonic
decreasing, hence it has a limit $c \geq 0$. Choose $n$ such that
$0 \leq \displaystyle\frac{c_{2^n}}{2^n} - c < \varepsilon$. For
any integer $p >\displaystyle \frac{2^nc_1}{\varepsilon}$ write
$$p = 2^{m_1}+ \cdots + 2^{m_s} + 2^{m_{s+1}} + \cdots +
2^{m_l},$$
 where $m_1 < \cdots < m_s < n \leq m_{s+1} < \cdots < m_l$. From the
 condition $c_{a+b} \le c_a + c_b$ we can deduce that
  $\displaystyle\frac{c_p}{p} <
 c+2\varepsilon$. Set $2^{m_l+2} = p+p'$. Since $c_p
 \geq c_{p+p'} - c_{p'}$ and $p' > \displaystyle \frac{2^nc_1}{\varepsilon}$,
  we have
  $\displaystyle\frac{c_p}{p} > c-8\varepsilon$. Hence the limit exists.
 Moreover, for each $n >0$, the sequence
 $\{\displaystyle\frac{c_{2^ln}}{2^ln}\}$ is
 decreasing and converges to $c$. Hence we get $c_n \geq cn$. \sk

\begin{Corollary} \label{a-filter} Let $(R,\mm)$ be a
local ring. Let $\{I_n\}$ be a filtration of ideals in $R$ with
$\dim R/I_n = d$ for a fixed integer $d$ and all $n \ge 0$. 
Then, for all $n \ge 0$, 
$$a(G(R/I_n)) \le [r(R/I_1)+1]n-d-1.$$ \end{Corollary}

\begin{pf} By Proposition \ref{bound-reduction} we have
$a(G(R/I_n))\le r(G(R/I_n))-d$. But $r(G(R/I_n)) = r(R/I_n)$ as
by Lemma \ref{comparison-reduction}. Hence the conclusion follows
from Theorem \ref{asymptotic-reduction}(i). \end{pf} \sk

Let $I$ be an arbitrary ideal in the local ring $R$. According to
Theorem \ref{asymptotic-reduction} and Corollary \ref{a-filter},
the reduction number $r(G(R/I^n))$ and the $a$-invariant
$a(G(R/I^n))$ are bounded by a linear function. On the other hand, 
each associated graded ring $G(R/I^n)$ can be represented as a quotient ring of the associated 
graded ring $G(R)$ by a homogeneous ideal, say $I^*_n$. Since $\{I^*_n\}$ is a filtration of 
ideals in $G(R)$, one may ask whether $\reg(G(R/I^n))$ is bounded by a linear function. We
can give a positive answer to this question in the case $\dim R/I
\le 1$. \sk

\begin{Theorem}\label{Reg3} Let $(R,\mm)$ be a local ring and
$I \subset R$ an ideal with $\dim R/I \leq 1$. Then there exists a
constant $c \ge 0$ such that for all $n>0$,
$$\reg(G(R/I^n))\leq cn.$$
\end{Theorem}

\begin{pf} The statement is trivial if $\dim R/I =0$. Let $\dim R/I =1$.
Write $I^n = J_n \cap Q_n,$ where $J_n$ is the
intersection of all primary components whose associated primes
are different from $\mm$ and $Q_n$ is the $\mm$-primary component
in an irredundant primary composition of $I^n$ ($Q_n = R$ if
$\mm$ is not an associated prime of $I^n$). By [Sw, Theorem 3.4],
there is a constant $a$ such that for all $n \ge 0$, there is an
irredundant decomposition of $I^n$ with $\mm^{an} \subseteq Q_n$.
>From this it follows that
$$I^n = J_n \cap (\mm^{an}+I^n).$$
Let $z \in {\frak m}$ be an element such that $z$ generates a minimal
reduction of $R/I^n$ with $r_{(z)}(R/I^n) = r(R/I^n)$. We have $J_n:z = J_n$ because $z$ is a 
non-zerodivisor on $J_n$.
Therefore, $I^n:z = J_n \cap [(\mm^{an}+I^n):z].$ Hence
$$(I^n:z) \cap (\mm^{an}+I^n) = J_n \cap (\mm^{an}+I^n) = I^n.$$
Applying Theorem \ref{Reg2} to $R/I^n$ with $\qq = (z,I^n)/I^n$ we
see that
$$\reg(G(R/I^n)) \le \max\{an,r(R/I^n)\}.$$
By Theorem \ref{asymptotic-reduction}(i), $r(R/I^n) \le
[r(R/I)+1]n-1$. Therefore, if we put $c = \max\{a,r(R/I)+1\}$,
then $\reg(G(R/I^n)) \le cn$ for all $n \ge 0$. \end{pf}\smallskip

\noindent{\it Remark}. If $d := \dim R/I < \dim R$, we have
$\displaystyle \lim_{n\rightarrow \infty} \frac{r(R/I^n)}{n} \ge
1$ and therefore $\displaystyle \lim_{n\rightarrow \infty}
\frac{\reg(G(R/I^n))}{n} \ge 1$ if it exists. Indeed, let $z_1,...,z_d \in \mm$
be elements such that $\underline{z}=\{ z_1,...,z_d\}$ generates a
minimal reduction of the maximal ideal of $R/I^n$. If ${\frak
m}^{n-1} \subseteq (\underline{z}){\frak m}^{n-2} + I^n$, then
${\frak m}^{n-1} = (\underline{z}){\frak m}^{n-2}$. Hence
$(\underline{z})$ is a minimal reduction of $\frak m$, a
contradiction to the assumption $d < \dim R$. So we have
$r_{\underline{z}}(R/I^n) \geq n-1$, which implies that
$\displaystyle \lim_{n\rightarrow \infty} \frac{r(R/I^n)}{n} \ge
1$. \sk

Now we will construct an example with $\reg(G(R/I^n)) = r(R/I^n)$ for all $n \ge 1$ such that 
$\displaystyle \lim_{n\rightarrow
\infty} \frac{r(R/I^n)}{n}$ is any given rational number not
less than $1$ (cf. [Cu]). As a consequence, $r(R/I^n)$ and
$\reg(G(R/I^n))$ need not to be linear functions for large $n$.
According to the result of [CHT], this implies that the Rees algebra
of the filtration $\{I^*_n\}$ in the associated graded ring $G(R)$ is not Noetherian in this 
case. \sk

\noindent{\bf Example.}  Let $R = k[[x,y,z_1,\ldots,z_d]]/(x^p-y^q)$, 
where $p < q$ are two different prime numbers, and $I = (x)$. 
Then $\dim R/I^n = d$. It is clear that $G(R/I^n)$ is a Cohen-Macaulay ring. Hence 
$\reg(G(R/I^n)) = r_\qq(R/I^n) = r(R/I^n)$ for any minimal reduction $\qq$ of the maximal 
ideal of $R/I^n$. We will
show that $\displaystyle \lim_{n\rightarrow \infty} \frac{r(R/I^n)}{n} =q/p.$ \par
Choose $\qq = (z_1,\ldots,z_d)$. Set $A = k[[t^p,t^q]]$. Then
$$r_\qq(R/I^n) = \max\{ s|\ (t^p,t^q)^sA \not\subseteq t^{qn}A\}.$$
For each integer $n$ write $n = pr + s$, where
$r \geq 0$ and $0\leq s\leq p-1$. Since $(t^p)^{qr + s} =
t^{pqr+ps} \not\in (t^q)^{pr+s}A$, we have $r_\qq(R/I^n) \geq qr +s$. 
On the other hand, let $l \geq qr + qs +p-1$. Assume that $u \in
(t^p,t^q)^lA$. Then $u$ is a sum of elements of the form $(t^p)^a(t^q)^bv$,
where $a+b \geq l$ and $v \in A$. Since $(t^q)^p = (t^p)^q$ and
$p < q$, we may assume that $b\leq p-1$ and $a \geq l-(p-1) \geq qr+qs$. Then
$$ \begin{array}{ll}
(t^p)^a(t^q)^bv & = (t^p)^{qr+qs} (t^p)^{a-qr-qs}(t^q)^bv\\
& = t^{pqr + qs}(t^q)^{ps-s}(t^p)^{a-qr-qs}(t^q)^bv\\
& \in (t^q)^{pr+s}A.
\end{array}$$
Hence $u \in t^nA$ and $r_\qq(R/I^n) < qr + qs +p-1$. Together with the inequality $qr+s \leq 
r_\qq(R/I^n)$, this implies that
$$\displaystyle \lim_{n\rightarrow \infty} \frac{r(R/I^n)}{n} =
\displaystyle \lim_{n\rightarrow \infty} \frac{r_\qq(R/I^n)}{n} =
\frac{q}{p}.$$ \par Using the argument above we can eventually compute
$r(R/I^n)$ for given $p, q$ and $n$. For instance, if $p=2$ and
$q=3$, then
 $$r(R/I^n) = \begin{cases} 3l \ \ & \text{if} \ n = 2l,\\
 3l + 2 \ \ & \text{if} \ n = 2l+1. \end{cases}$$\sk

\section{Asymptotic regularity of symbolic powers}\sk

Let $R = k[x_1,\ldots,x_n]$ be a polynomial ring over a field $k$
and $\mm$ the maximal graded ideal of $R$. Let $I$ be an arbitrary
homogeneous ideal of $R$. Given an ideal (or a subset) $J \subset
R$, we set
$$I:\langle J \rangle = \cup_{n\ge 0}I:J^n.$$ It is easily seen
that $I:\langle J \rangle$ is the intersection of all components
in a primary decomposition of $I$ whose associated primes do not
contain $J$. \sk

\noindent{\bf Definition.} For every integer $n \ge 1$ we call
the ideal
$$I^{(n)} := I^n:\langle J \rangle.$$
the $n$th {\it $J$-symbolic power} of $I$.\sk

This notion generalizes the (ordinary) symbolic power and the
saturated power of $I$. Let $A^*(I)$ denote the union of the sets
of the associated primes of $I^n$ for all $n \ge 0$. It is
well-known that $A^*(I)$ is a finite set (see e.g. [Mc]) and that
$A^*(I) \supseteq \Min(I)$, where $\Min(I)$ denotes the set of
all minimal associated primes of $I$. If $J$ is the intersection
of all primes of $A^*(I)\setminus\Min(I)$, then $I^{(n)}$ is the
(ordinary) $n$th symbolic power of $I$. If $J = \mm$, then
$I^{(n)}$ is the saturation of $I^n$. \sk

\noindent{\it Remark.} The symbolic Rees algebra $\oplus_{n \ge
0}I^{(n)}t^n$ needs not to be a Noetherian ring. For instance, if $I$ is the defining ideal of 
a set of $2r+1$ points on a rational curve in ${\Bbb P}^r$, $r \ge 2$, then 
$$\reg((I^n)^\sat) = n+1+ \left[\frac{n-2}{r}\right]$$ is not a linear function [CTV, 
Proposition 7]. Thus, the saturated Rees algebra $\oplus_{n\ge 0}(I^n)^\sat t^n$ is not 
Noetherian, according to [CHT, Lemma 3.3 and Theorem 3.4]. \sk

Since $\{I^{(n)}\}$ is a filtration of homogeneous ideals with $\dim R/I^{(n)} = \dim 
R/I^{(1)}$, the $a$-invariant of $R/I^{(n)}$ is always bounded by a linear function of $n$. 
\sk

\begin{Proposition} \label{a-symbolic} Let $I$ be an arbitrary
homogeneous ideal of $R$ and $d = \dim R/I^{(1)}$. For all $n \ge
0$ we have
$$a(R/I^{(n)}) \le [r(R/I^{(1)})+1]n - d -1.$$ \end{Proposition}

\begin{pf} By Theorem \ref{asymptotic-reduction}(i) we have 
$$r(R/I^{(n)}) \le [r(R/I^{(1)})+1]n - 1.$$
Hence the conclusion follows from Proposition
\ref{bound-reduction}. \end{pf} \sk

Inspired of Theorem \ref{ordinary} and Theorem \ref{saturation} one may ask whether there is a 
linear upper bound for $\reg(I^{(n)})$. We shall see that this question has a positive answer 
if $\dim R/J \le 1$ or if $I$ is a
monomial ideal. \sk

\begin{Proposition} \label{comparison} Let $J \subset R$ be an ideal
with $\dim R/J \le 1$. For all $n \ge 0$ we have
$$\reg(I^{(n)}) \le \reg(I^n)^\sat \le \reg(I^n).$$
\end{Proposition}

\begin{pf} Since $I^{(n)} = \cup_{r\ge
0}I^n:J^r$, the quotient module $I^{(n)}/I^n$ is annihilated by
some power of $J$. It follows that $\dim I^{(n)}/I^n \le \dim R/J
\le 1$. Therefore, $H_{\mm}^i(I^{(n)}/I^n) = 0$ for all $i \ge 2$.
Now, from the exact sequence
$$0 \To I^{(n)}/I^n \To R/I^n \To R/I^{(n)} \To 0$$
we can deduce that there is a surjective map $H_{\mm}^1(R/I^n) \To
H_{\mm}^1(R/I^{(n)})$ and $H_{\mm}^i(R/I^{(n)}) \cong
H_{\mm}^i(R/I^n)$ for $i \ge 2$. It follows that $a_1(R/I^{(n)})
\le a_1(R/I^n)$ and $a_i(R/I^{(n)}) = a_i(R/I^n)$ for $i \ge 2$.
Since $\mm$ is not an associated prime of $I^{(n)}$, we have
$H_{\mm}^0(R/I^{(n)}) = 0$ and therefore $a_0(R/I^{(n)}) =
-\infty$. So we can conclude that
$$\reg(R/I^{(n)}) = \max\{a_i(R/I^{(n)})+i|\ i > 0\}
\le \max\{a_i(R/I^n)+i|\ i > 0\}.$$ On the other hand, since
$(I^n)^\sat/I^n$ is a module of finite length,
$H_{\mm}^i(R/(I^n)^\sat) = H_{\mm}^i(R/I^n)$ for $i
> 0$. Since $H_{\mm}^0(R/(I^n)^\sat) = 0$, we have
$$\reg(R/(I^n)^\sat) = \max\{a_i(R/I^n)+i|\ i > 0\}.$$
So we obtain $\reg(R/I^{(n)}) \le \reg(R/(I^n)^\sat) \le
\reg(R/I^n)$, which implies the conclusion.
\end{pf} \sk

\begin{Theorem} \label{symbolic-general} Let $J \subset R$ be an ideal
with $\dim R/J \le 1$. Then there exists a constant $e$ such that
for all $n > 0$,
$$\reg(I^{(n)}) \le d(I)n+e.$$ \end{Theorem}

\begin{pf} This follows from Theorem \ref{ordinary}
and Proposition \ref{comparison}. \end{pf} \sk

\noindent{\it Remark.}  We may apply Theorem \ref{saturation} to
obtain a better bound for $\reg(I^{(n)})$. \sk

The following example shows that if $\dim R/J \ge 2$, the inequality $\reg(I^{(n)}) \le 
\reg(I^n)$ does not hold even for monomial ideals. \sk

\noindent{\bf Example.}  Let $R = k[x_1,x_2,x_3,x_4]$ and
\begin{eqnarray*} I & = & (x_1,x_2^2) \cap (x_3,x_4^2) \cap (x_1,x_3)\\ 
& = &(x_1x_3,x_1x_4^2,x_2^2x_3). \end{eqnarray*} Let $J = (x_1,x_3)$. Then $\dim R/J = 2$ and
$$I^{(n)} = (x_1,x_2^2)^n \cap (x_3,x_4^2)^n$$ for all $n \ge 1$. 
Applying the graded version of Theorem \ref{Reg2} to the rings $R/I^n$ and $R/I^{(n)}$ with 
$z_1 = x_2+x_4, z_2 = x_1+x_2+x_3$ we can show that $\reg(R/I^n) = 3n-1$ and $\reg(R/I^{(n)}) 
= 4n-1$ for all $n \ge 1$. From this it follows that $$\reg(I^n) = 3n < 4n = \reg(I^{(n)}).$$ 

Now we will apply Theorem \ref{symbolic-general} to the ordinary
symbolic powers. \sk

\begin{Corollary} \label{symbolic=2} {\rm (cf. [Ch, Corollary 7]) }
Let $\dim R/I \le 2$. Then, for the ordinary symbolic powers of
$I$, there exists a constant $e$ such that for all $n \ge 0$,
$$\reg(I^{(n)}) \le d(I)n+e.$$ \end{Corollary}

\begin{pf} For the ordinary symbolic powers, we can choose $J$ to
be the intersection of the primes of $A^*(I)\setminus \Min(I)$.
>From this it follows that $\dim R/J \le \dim R/I - 1 \le 1$.
Hence the conclusion follows from Theorem \ref{symbolic-general}.
\end{pf} \sk

\noindent{\it Remark.} Chandler proved that if $\dim R/I \le 2$, then $\reg(I^{(n)}) \le 
\reg(I)n$ for all $n \ge 1$.\sk

\begin{Corollary} \label{symbolic-singular} Assume that the
singular locus of $R/I$ has dimension $\le 1$. Then, for the
ordinary symbolic powers of $I^n$, there exists a constant $e$
such that for all $n \ge 0$,
$$\reg(I^{(n)}) \le d(I)n+e.$$ \end{Corollary}

\begin{pf} Since the singular locus of $R/I$ has dimension $\le
1$, $(R/I)_\pp $ is a regular local ring for all primes $\pp
\supseteq I$ with $\dim R/\pp  > 1$. For any such $\pp$ it follows
that $I_\pp $ is a complete intersection and hence $I^n_\pp $ is
an unmixed ideal. As a consequence, $I^n$ has no embedded
associated primes $\frak p$ with $\dim R/\pp  > 1$. That means
$\dim R/\pp < 1$ for all primes $\pp  \in A^*(I)\setminus
\Min(I)$. If we choose $J$ to be the intersection of the primes
of $A^*(I) \setminus \Min(I)$, then $\dim R/J \le 1$. Therefore,
the conclusion follows from Theorem \ref{symbolic-general}.
\end{pf} \sk

In the case $\dim R/J = 2$ we can only show that if $d(I^{(n)})$, the maximal degree of the 
defining equations of $I^{(n)}$, is bounded by a linear function, then $\reg(I^{(n)})$ is 
bounded by a linear function, too.
Note that if $\reg(I^{(n)})$ is bounded by a linear function, then so is 
$d(I^{(n)})$ and that the latter condition can be checked often. \sk

\begin{Proposition} \label{symbolic-generator} Let $J \subset R$ be an ideal with $\dim R/J = 
2$. Assume that there is an integer $a$ such that
$d(I^{(n)}) \le an$ for all $n\ge 1$. Then there is an integer $c$ such that  for all $n \ge 
1$, 
$$\reg(I^{(n)}) \le cn.$$ \end{Proposition}
 
\begin{pf} From the proof of Proposition \ref{comparison} we see that $a_0(R/I^{(n)}) = 
-\infty$ and $$a_1(R/I^{(n)}) \le \max\{a_1(R/I^n),a_2(I^{(n)}/I^n)\}.$$ Moreover, since $\dim 
I^{(n)}/I^n \le \dim R/J = 2$, we have $a_i(I^{(n)}/I^n) = 0$ for $i \ge 3$. >From this it 
follows that $a_2(R/I^{(n)}) \le a_2(R/I^n)$ and $a_i(R/I^{(n)}) \le a_i(R/I^n)$ for $i \ge 
3$. Thus, it suffices to show that $a_2(I^{(n)}/I^n)$ is bounded by a linear function. \par  
Write $I^n = I^{(n)}\cap Q_n$, where $Q_n$ is the intersection of all components of a primary 
decomposition of $I^n$ whose associated primes contain $J$. Then $$I^{(n)}/I^n \cong 
(I^{(n)}+Q_n)/Q_n.$$
By [Sw, Theorem 3.4] we may assume that $J^{bn} \subseteq Q_n$ for a fixed integer $b \ge 1$. 
Choose $z_1,z_2 \in R_1$ such that $\qq = (z_1,z_2)$ is a minimal reduction in $R/J^{bn}$ with 
$r_\qq(R/J^{bn}) = r(R/J^{bn})$. Then $r_\qq(R/J^{bn}) \le [r(R/J)+1]bn-1$ by Theorem 
\ref{asymptotic-reduction}.
Put $m = [r(R/J)+1]bn$. Then $R_m = (z_1,z_2,I^n)_m$. Since $d(I^{(n)}) \le an$, we have 
$$(I^{(n)})_{an+m} = R_m(I^{(n)})_{an} = (z_1,z_2,I^n)_m(I^{(n)})_{an}.$$
>From this it follows that $(I^{(n)}/I^n)_{an+m} = \big[(z_1,z_2)(I^{(n)}/I^n)\big]_{an+m}$. 
Hence $\qq$ is also a minimal reduction of $I^{(n)}/I^n$ with $r_\qq(I^{(n)}/I^n) \le an+m-1$. 
Thus, $r(I^{(n)}/I^n) \le an+m-1$. By Proposition \ref{bound-reduction}, this implies that 
$a_2(I^{(n)}/I^n) \le an+m-3$. Hence $a_2(I^{(n)}/I^n)$ is bounded by a linear function, as 
required. \end{pf} \sk

\begin{Corollary} Let $\dim R/I = 3$. Let $I^{(n)}$ be the ordinary $n$th symbolic power of 
$I$. Assume that there is an integer $a$ such that $d(I^{(n)}) \le an$ for all $n \ge 1$. Then 
there is an integer $c$ such that  for all $n \ge 1$, 
$$\reg(I^{(n)}) \le cn.$$ \end{Corollary}

\begin{pf} For the ordinary symbolic powers, we can choose $J$ to
be the intersection of the primes of $A^*(I)\setminus \Min(I)$. Hence
the assumption $\dim R/I = 3$ implies $\dim R/J \le 2$. \end{pf} \sk

In the following we will show that $\reg(I^{(n)})$ is bounded by a
linear function for all monomial ideals $I$ and arbitrary ideals
$J \supseteq I$. The proof is based on a result of Bruns and
Herzog which says that the regularity of a monomimal ideal $I$ is
bounded by the least common multiple of the generating monomials
of $I$ [BH2, Theorem 3.1]. We shall need the following simple
observation. \sk

\begin{Lemma} \label{divisor} Let $I$ be an arbitrary monomial
ideal. Let $f$ and $f'$ be the least common multiples of the
generating monomials of $I$ and $I:\langle J\rangle$,
respectively. Then $f'$ is a divisor of $f$. \end{Lemma}

\begin{pf} Let $f = x_1^{a_1}\cdots x_m^{a_m}$. Then $x_i^{a_i}$ is the
largest power of $x_i$ which appears in a generating monomial of
$I$, $i = 1,\ldots,m$. Using the formula
$$(I',g_1g_2) = (I',g_1) \cap (I',g_2),$$
where $I'$ is a  monomial ideal and $g_1$ and $g_2$ are monomials
which have no common divisor, we can find a primary decomposition
of $I$ such that $x_i^{a_i}$ is the largest exponent of $x_i$
which occurs as a generator of some primary component. The ideal
$I:\langle J\rangle$ is obtained from this decomposition by
deleting those primary components whose associated primes contain
$J$. Thus, the exponent of $x_i$ in every generating monomial of
$I:\langle J\rangle$ must be less than or equal to $a_i$. Hence $f$ is
a common multiple of the generating monomials of $I:\langle
J\rangle$.
\end{pf} \sk

\begin{Theorem} \label{symbolic-monomial} Let $I$ and $J$ be arbitrary
monomial ideals. Let $f$ be the least common multiple of the
generating monomials of $I$. Then, for all $n \ge 0$,
$$\reg(I^{(n)}) \le (\deg f)n - \height I + 1.$$ \end{Theorem}

\begin{pf} We first note that $I^{(n)}$ is the intersection of
some primary components of $I^n$. Let $g_n$ denote the least
common multiple of the generating monomials of $I^{(n)}$. It has
been shown in [HT, Theorem 3.1] (cf. [BH2, Theorem 3.1] for a
weaker form) that
$$\reg(I^{(n)}) \le \deg g_n -
\height I^{(n)} + 1 \le \deg g_n - \height I +1.$$ If we denote by
$f_n$ the least common multiple of the generators of $I^n$, then
$\deg g_n \le \deg f_n$ by Lemma \ref{divisor}. On the other hand,
since $f^n$ is a common multiple of the generating monomials of
$I^n$, $f_n$ is a divisor of $f^n$. Therefore, $\deg f_n \le \deg
f^n = (\deg f)n$. Hence $\reg(I^{(n)}) \le (\deg f)n - \height I
+ 1.$ \end{pf} \sk

\noindent{\it Remark.} The previous example shows that we do not always have the inequality 
$\reg(I:\langle J\rangle) \le \reg(I)$ for  monomial ideals $I$ and $J$. \sk

\section{Asymptotic regularity of initial ideals}\sk

Let $R=k[x_1,\ldots,x_m]$ be a polynomial ring over a field $k$
and $\mm$ the maximal graded ideal of $R$. Let $I$ be an
arbitrary homogeneous ideal of $R$. Let $\ini(I)$ denote the
initial ideal of $I$ with respect to a term order of $R$. The aim
of this section is to study the regularity of the initial ideal
$\ini(I^n)$. We will do it in a more general setting by
investigating the initial ideals of a weight order. \sk

Given a linear map $\lambda: {\Bbb Z}^m\to {\Bbb Z}$, we can
define a weight order on the monomials of $R$. Let
$\ini_{\lambda}(I)$ denote the initial ideal of $I$ with respect
to this monomial order. It can be shown that $\ini(I) =
\ini_{\lambda}(I)$ for a suitable choice of the integral weight function
$\lambda$ which depends on $I$ (see e.g. [E, p.~327]). The ideal
$\ini_{\lambda}(I)$ can be obtained as follows. \sk

Let $R[t]$ be a polynomial ring over $R$ in one variable $t$. For
any $g\in R[t]$, $g=\sum_ia_iu_i$, where the $u_i$ are monomials
and $0\neq a_i\in k$, we set $b(g)=\max\lambda(u_i)$ and define
$g^* \in R[t]$ as follows:
\[ g^* :=
t^{b(g)}g(t^{-\lambda(x_1)}x_1,\ldots,t^{-\lambda(x_r)}x_m).
\]
We denote by $I^*$ the ideal of $R[t]$ generated by $\{g^*|\ g\in
I\}$. It is known that $t$ is a non-zerodivisor modulo $I^*$ and that
$R[t]/(I^*,t) \cong R/\ini_{\lambda}(I)$. \sk

In order to study the regularity of $\ini_{\lambda}(I)$ we
consider the extension of $I^*$ in the localization $S = R
\otimes k[t]_{(t)}$. We will view $S$ as a standard graded
algebra over the local ring $k[t]_{(t)}$ ($\deg t = 0$). Let
$\widetilde{I} = I^*S$. Then $t$ is still a non-zerodivisor modulo
$\widetilde{I}$ and
$$S/(\widetilde{I},t) \cong
R/\ini_{\lambda}(I).$$ Since $\widetilde{I}$ and
$\ini_{\lambda}(I)$ are homogeneous ideals with respect to this
grading, we can define the regularity and the reduction number
for these ideals as in Section 1. \sk

\begin{Proposition} \label{comparison-ini}
With the above notations we have $$\reg(\ini_{\lambda}(I)) =
\reg(\widetilde{I}).$$
\end{Proposition}

\begin{pf} From the exact sequence $0 \To S/\widetilde{I} \overset t \To
S/\widetilde{I} \To R/\ini_{\lambda}(I) \To 0$ we get the long
exact sequence
$$ \cdots \To H_{S_+}^i(S/\widetilde{I}) \overset t
\To H_{S_+}^i(S/\widetilde{I}) \To H_{\mm}^i(R/\ini_{\lambda}(I))
\To H_{S_+}^{i+1}(S/\widetilde{I}) \To \cdots, $$ where the
homomorphisms are all of degree $0$. It follows that
$a_i(R/\ini_{\lambda}(I)) \le
\max\{a_i(S/\widetilde{I}),a_{i+1}(S/\widetilde{I})\}$. Hence
\begin{eqnarray*} \reg(\ini_{\lambda}(I)) & = &
\max\{a_i(R/\ini_{\lambda}(I))+i|\ i \ge 0\}\\ & \le &
\max\{a_i(S/\widetilde{I})+i|\ i \ge 0\}\ =\
\reg(S/\widetilde{I}),
\end{eqnarray*} which implies $\reg(\ini_{\lambda}(I)) \le
\reg(\widetilde{I})$. \par To prove the converse, let $a_i =
a_i(R/\ini_{\lambda}(I))$. Then $H_{\mm}^i(R/\ini_{\lambda}(I))_n
= 0$ for $n > a_i$. Therefore, $tH_{S_+}^i(S/\widetilde{I})_n =
H_{S_+}^i(S/\widetilde{I})_n$ for $n > a_i$. Since
$H_{S_+}^i(S/\widetilde{I})_n$ is a finitely generated module
over the local ring $k[t]_{(t)}$, this implies
$H_{S_+}^i(S/\widetilde{I})_n = 0$ by Nakayama's lemma. Thus,
$a_i(S/\widetilde{I}) \le a_i$. Hence
$$\reg(S/\widetilde{I}) \le \max\{a_i+i|\ i \ge 0\} =
\reg(R/\ini_{\lambda}(I)),$$ which implies $\reg(\widetilde{I})
\le \reg(\ini_{\lambda}(I)).$  \end{pf} \sk

\begin{Lemma} \label{reduction-ini}
$r(R/\ini_{\lambda}(I)) = r(S/\widetilde I)$.
\end{Lemma}

\begin{pf} First we note that $R/\ini_{\lambda}(I) =
S/(\widetilde I,t)$ is the fiber ring of $S/\widetilde I$.
Therefore, a reduction $J$ of $S/\widetilde I$ is minimal if and
only if $J$ is a minimal reduction of $R/\ini_{\lambda}(I)$.
Moreover, we have $r_J(S/\widetilde I) =
r_J(R/\ini_{\lambda}(I))$. From this it follows that
$r(R/\ini_{\lambda}(I)) = r(S/\widetilde I)$.
\end{pf} \sk

According to Proposition \ref{comparison-ini}, to estimate
$\reg(\ini_{\lambda}(I))$ we only need to estimate
$\reg(\tilde{I})$. To compute the local cohomology modules of
$S/\tilde{I}$ we shall need the following observation. \sk

\begin{Lemma} \label{a0} $a_0(S/\tilde{I}) \le a_0(R/I)$. \end{Lemma}

\begin{pf} We have $H_{\mm}^0(R/I) = I:\langle
\mm\rangle/I$ and, since $S_+$ is generated by $\mm$, $
H_{S_+}^0(S/\widetilde{I}) = \widetilde{I}:\langle \mm
\rangle/\widetilde{I}$. So we have to prove that
$[\widetilde{I}:\langle \mm \rangle/\widetilde{I}]_n = 0$ if
$[I:\langle \mm\rangle/I]_n = 0$. \par Let $T = R[t,t^{-1}]$. Then
$T$ can be viewed as a standard graded algebra over
$k[t,t^{-1}]$. Let $\phi$ be the graded automorphism of $T$
determined by $\phi(x_i) = t^{-\lambda(x_i)}x_i$, $i =
1,\ldots,m$. Let $\phi(I)$ denote the ideal of $R[t,t^{-1}]$
generated by the elements $\phi(g)$, $g \in I$. It is easy to
check that
$$I^* = \phi(I) \cap R[t].$$
>From this it follows that
$$I^*:\langle \mm\rangle = (\phi(I):\langle \mm \rangle) \cap R[t]
= \phi(I:\langle \mm \rangle) \cap R[t].$$ If we consider $R[t]$
as a standard graded algebra over $k[t]$, then these formulae also
hold for the corresponding graded pieces. Thus, if $(I:\langle
\mm\rangle)_n = I_n$, then $(I^*:\langle \mm\rangle)_n =
(I^*)_n$; hence $(\widetilde{I}:\langle \mm \rangle)_n =
\widetilde{I}_n$. \end{pf} \sk

Now, we will use the above relationships between
$\ini_{\lambda}(I^n)$ and $\widetilde{I^n}$ to study the
asymptotic regularity of $\ini_{\lambda}(I^n)$.\sk

\begin{Theorem} \label{asymptotic-ini} Let $\dim R/I \le 1$.
Then\par {\rm (i) } $\reg(\ini_{\lambda}(I^n)) \le
\reg(\ini_{\lambda}(I))n$ for all $n \ge 0$,
\par {\rm (ii) } $c = \displaystyle
\lim_{n\to\infty}\frac{\reg(\ini_{\lambda}(I^n))}{n}$ exists with
$c \le \max\{d(I),r(R/\ini_{\lambda}(I))+1\}$.
\end{Theorem}

\begin{pf} By Proposition \ref{comparison-ini} and Lemma
\ref{reduction-ini} we may replace $\ini_{\lambda}(I^n)$ by
$\widetilde{I^n}$ in the above statements. Note that
$S/(\widetilde{I^n},t) \cong R/\ini_{\lambda}(I^n)$ is the fiber
ring of $S/\widetilde{I^n}$. Since $\dim R/\ini_{\lambda}(I^n) =
\dim S/I^n = \dim S/I \le 1$, $a_i(S/\widetilde{I^n}) = 0$ for $i
> 2$. Therefore, using Proposition \ref{bound-reduction} we
obtain the formula
$$\reg(S/\widetilde{I^n}) =
\max\{a_0(S/\widetilde{I^n}),r(S/\widetilde{I^n})\}.$$ By Lemma
\ref{a0}(ii) we have $a_0(S/\widetilde{I^n}) \le a_0(R/I^n) \le
\reg(R/I^n).$ On the other hand, $\reg(I^n) \le \reg(I)n$ by
[Ch] and [GGP]. From this it follows that
$\reg(R/I^n) \le [\reg(R/I)+1]n-1$. Since $R/\ini_{\lambda}(I)$
is a flat deformation of $R/I$ [E, Theorem 15.17], $\reg(R/I) \le
\reg(R/\ini_{\lambda}(I)) = \reg(S/\tilde I)$. Therefore,
$$a_0(S/\widetilde{I^n}) \le [\reg(S/\tilde I)+1]n-1.$$
Since $\{\widetilde{I^n}\}$ is a filtration of ideals with $\dim
S/(\widetilde{I^n},t) = \dim R/I$ for all $n \ge 0$, we can apply
the graded version of Theorem \ref{asymptotic-reduction}(i) and
obtain $r(S/\widetilde{I^n}) \le [r(S/\tilde{I})+1]n-1$. On the
other hand, $r(S/\tilde{I}) \le \reg(S/\tilde I)$ by Proposition
\ref{bound-reduction}. Therefore, we also have
$$r(S/\widetilde{I^n}) \le [\reg(S/\tilde I)+1]n-1.$$
Summing up, we get $\reg(S/\widetilde{I^n}) \le
[\reg(S/\tilde I)+1]n-1,$ which implies (i).\par By the above
arguments we have
$$\reg(S/\widetilde{I^n}) \le
\max\{\reg(R/I^n),r(S/\widetilde{I^n})\} \le
\reg(S/\widetilde{I^n}).$$ This implies $\reg(S/\widetilde{I^n}) =
\max\{\reg(R/I^n),r(S/\widetilde{I^n})\}$. By Theorem
\ref{ordinary} we have  $\displaystyle
\lim_{n\to\infty}\frac{\reg(R/I^n)}{n} \le d(I)$. By Theorem
\ref{asymptotic-reduction}, $c_1 = \displaystyle
\lim_{n\to\infty}\frac{r(S/\widetilde{I^n})}{n}$ exists with $c_1
\le r(S/\widetilde I)+1$. So we get
\begin{eqnarray*} \lim_{n\to\infty}\frac{\reg(S/\widetilde{I^n})}{n}
& = & \max\big\{ \lim_{n\to\infty}\frac{\reg(R/I^n)}{n},
\lim_{n\to\infty}\frac{r(S/\widetilde{I^n})}{n}\big\}\\ & \le &
\max\{d(I),r(S/\widetilde I))+1\},\end{eqnarray*} which implies
(ii). The proof of Theorem \ref{asymptotic-ini} is now complete.
\end{pf} \sk

We can not simply translate Theorem \ref{asymptotic-ini} for
initial ideals (with respect to a term order). For, we can not always find an integral weight 
function $\lambda$ such that $\ini(I^n) =
\ini_{\lambda}(I^n)$ for all $n \ge 0$.\sk

\begin{Theorem} \label{ini}
Let $I$ be a homogeneous ideal with $\dim R/I \le 1$.
Then \par {\rm (i) } $\reg(\ini(I^n)) \le \reg(\ini(I))n$ for all
$n \ge 0$, \par {\rm (ii) } $c = \displaystyle
\lim_{n\to\infty}\frac{\reg(\ini(I^n))}{n}$ exists with $c \le
\max\{d(I),r(R/\ini(I))+1\}$.
\end{Theorem}

\begin{pf} For every $n \ge 0$ we can find an integral weight function $\lambda$ such that 
$\ini(I^r) = \ini_{\lambda}(I^r)$ for all $r \le n$.
Therefore, (i) follows from Theorem \ref{asymptotic-ini}(i). The
proof for (ii) is similar to the proof for Theorem
\ref{asymptotic-ini}(ii). Indeed, we have
$$\reg(S/\widetilde{I^n}) = \max\{\reg(R/I^n),r(S/\widetilde{I^n})\}.$$
By Proposition \ref{comparison-ini} and Lemma \ref{reduction-ini},
this implies
$$\reg(R/\ini(I^n)) =
\max\{\reg(R/I^n),r(R/\ini(I^n))\}.$$ Hence the conclusion follows
from Theorem \ref{ordinary} and Theorem
\ref{asymptotic-reduction}. \end{pf} \sk

\noindent{\it Remark.} If $\dim R/I \ge 2$, Sturmfels [St] gave
an example showing that the inequality $\reg(I^n) \le
\reg(I)n$ does not hold for monomial ideals. Therefore, Theorem
\ref{asymptotic-ini}(i) and Theorem \ref{ini} (i) do not hold in 
this case. \sk

We do not know whether $\reg(\ini(I^n))$ is bounded by a linear
function in general. By [BH2, Theorem 3.1] and [HT, Theorem 3.1], this is the case if and only 
if $d(\ini(I^n))$ is bounded by a linear function. 
Apart from the above positive answer for the
case $\dim R/I \le 1$, we can only show that the $a$-invariant of
$R/\ini(I^n)$ is bounded by a linear function. \sk

\begin{Proposition} \label{a-ini} Let $I$ be an arbitrary homogeneous
ideal of $R$ and $d = \dim R/I$. For all $n \ge 0$ we have
$$a(R/\ini(I^n)) \le [r(R/\ini(I))+1]n - d - 1.$$
\end{Proposition}

\begin{pf} Since $\{\ini(I^n)\}$ is a filtration of
homogeneous ideals with $\dim R/\ini(I^n) = \dim R/I$, we can
apply Theorem \ref{asymptotic-reduction}(i) and obtain
$$r(R/\ini(I^n)) \le [r(R/\ini(I))+1]n - 1.$$
Hence the conclusion follows from Proposition
\ref{bound-reduction}. \end{pf} \sk

We shall see that the ideal $\widetilde{I^n}$ may be considered
as a symbolic power of $\widetilde{I}$. \sk

\begin{Lemma} \label{symbolic-ini} $\widetilde{I^n} =
\widetilde{I}^n:\langle t \rangle$ for all $ n \ge 0$
\end{Lemma}

\begin{pf} It suffices to show that $(I^n)^* = (I^*)^n:\langle t
\rangle$. Since $t$ is a non-zerodivisor on $(I^n)^*$ and since
$(I^n)^* \supseteq (I^*)^n$, we have $(I^n)^* = (I^n)^*:\langle
t\rangle \supseteq (I^*)^n:\langle t\rangle.$  Conversely, if $f
= f_1\cdots f_n \in I^n$ with $f_1,\ldots,f_n \in I$, then we can
find an integer $a$ such that $t^af^* = f_1^*\cdots f_n^* \in
(I^*)^n$. From this it follows that for any element $g \in I^n$,
we have $t^ag^* \in (I^*)^n$ for all $a$ large enough, hence
$g^* \in (I^*)^n:\langle t\rangle$. So we can conclude that
$(I^n)^* = (I^*)^n:\langle t \rangle$. \end{pf} \sk

The above observation has the following interesting consequence.
\sk

\begin{Proposition} Let $I\subset R$ be a homogeneous ideal of
codimension $s$. Then for all $n \ge 0$ we have
$$\ini(I^{sn})\subseteq \ini(I)^n.$$
\end{Proposition}

\begin{pf} For a fixed natural number $n$ there exists an integral weight
function $\lambda$ such that $\ini(I)=\ini_\lambda(I)$ and
$\ini(I^{sn})=\ini_\lambda(I^{sn})$. By Lemma \ref{symbolic-ini},
$\widetilde{I^{sn}}= \tilde{I}^{sn}: \langle t\rangle$. Since $t$ is a non-zerodivisor modulo 
$\widetilde I$, we have $\tilde I^{sn}:\langle t\rangle \subseteq \tilde I^{(sn)}$. By a 
recent result of Ein,
Lazardsfeld and Smith [ELS] (see also [HH]), for any ideal $J\subset S$ of codimension $s$ and 
all $n \ge 0$ we have $\tilde J^{(sn)}\subset \tilde J^n$ for the ordinary symbolic powers. 
Therefore we get
\[ \widetilde{I^{sn}}= \tilde{I}^{sn}: \langle t\rangle \subseteq
\tilde{I}^{(sn)}\subseteq \tilde{I}^n. \] 
This implies $(\widetilde{I^{sn}},t) \subseteq (\tilde I^n,t)$.
Hence $\ini(I^{sn})\subseteq \ini(I)^n$. \end{pf} \sk

Finally, we will give some examples which show that the initial
Rees algebra $\oplus_{n\ge 0}\ini(I^n)t^n$ needs not to be
Noetherian.\sk

\noindent {\bf Example.} Let $k$ be any field, $I\subset k[x,y]$
the ideal generated by $x^3+y^3$ and $xy^2$, and let $<$ be any
term order with $x>y$. The generators of $\ini(I^n)$ in degree
$3n$ are $x^{n+2i}y^{2n-2i}$ for $i=0,\ldots, n$. Considering
$S$-pairs of the generators of $I^n$, we see that
$\ini(I^n)$ has no generators in degree $3n+1$. But on the other hand $y^{3n+2}$ is a minimal
generator of $\ini(I^n)$. In fact, $y^5=y^2(x^3+y^3)-x^2(xy^2)$ is an element of $I$, and
$y^{3n+2}=y^{3n-1}(x^3+y^3)-xy^{3n-5}(xy^2)^2$ belongs to $I^n$, since by the induction 
hypothesis we may assume  that
$y^{3n-1}\in I^{n-1}$ and $y^{3n-5}\in I^{n-2}$. Now it is easily seen that $y^{3n+2}t^n$ 
belongs to a minimal set of algebra generators of the Rees algebra $\oplus_{n\ge 
0}\ini(I^n)t^n$.\sk

The next example shows that not even $\oplus_{n\ge 0}\gin(I^n)t^n$ need to be finitely 
generated. \sk

\noindent {\bf Example.}  Let $k$ be a field of  characteristic 0
and $R = k[x_1,...,x_m]$. Let $I=(f_1,...,f_m)$ be an ideal
generated by a regular sequence of forms of the same degree $a>1$
in $R$. Since $\dim R/I = 0$, we have $\reg(R/I^n) = a(R/I^n) = (n-1)a +
m(a-1)$ for all $n>0$. Fix a term oder. By extending the field
$k$ to $k(u)$, where $u$ is a set of new variables, we may assume
that $\ini(I^n) = \gin(I^n)$ for all $n$. Each $\gin(I^n)$ has a
minimal generator of the form $x_m^{c_n}$. Since $\gin(I^n)$ is a
stable ideal and since $\dim R/I = 0$, we have
$$c_n = \reg(\gin (I^n)) = \reg(I^n) = (n-1)a + m(a-1)+1.$$
Were $\oplus_{n\geq 0} \gin(I^n)t^n$ finitely generated,  there
would exist a number $n_0$ such that $$\gin(I^n) =
\sum_{n_1+\cdots + n_s = n,\ n_1,...,n_s \leq n_0}
\gin(I^{n_1})\cdots \gin(I^{n_s})$$
 for all $n>0$. From this it follows that
  $$c_n = \min\{c_{n_1} + \cdots + c_{n_s}|
  \ n_1+\cdots + n_s = n, \ n_1,...,n_s \leq n_0 \}.$$
 We have $c_{n_1} + \cdots + c_{n_s} = (n-1)a + s[m(a-1) + 1] $.
 For $n > n_0$, we must have $s \ge 2$ for all sequences
 $n_1,\ldots,n_s \le n_0$ with $n_1+\cdots + n_s =
 n$. So we obtain $c_{n_1} + \cdots + c_{n_s} \geq (n-1)a + 2[m(a-1)+ 1]
  > c_n,$ a contradiction. \smallskip

\section*{References}\sk

\noindent [BEL] A.~Bertram, L. Ein and R. Lazarsfeld, Vanishing
theorems, a theorem of Severi, and the equations defining
projective varieties, J. Amer. Math. Soc. 4 (1991), no. 3,
587--602. \par

\noindent [BH1] W. Bruns and J. Herzog, Cohen-Macaulay rings,
Cambridge, 1998.\par

\noindent [BH2] W. Bruns and J. Herzog, On multigraded
resolutions, Math. Proc. Cambridge Phil. Soc. 118 (1995),
245-275.\par

\noindent [CTV] M.~Catalisano, N.V.~Trung and G. Valla, A sharp bound
for the regularity index of fat points in general position,
Proc. Amer. Math. Soc. 118 (1993), 717-724. \par

\noindent [Ch] K.~A. Chandler, Regularity of the powers of an
ideal, Commun. Algebra 25 (1997), 3773-3776.\par

\noindent [Cu] S.~D. Cutkosky, Irrational asymptotic behaviour of
Castelnuovo-Mumford regularity, J. Reine Angew. Math. 522 (2000),
93--103. \par

\noindent [CEL] S.~D. Cutkosky, L. Ein and R. Lazarsfeld,
Positivity and complexity of ideal sheaves, Preprint. \par

\noindent [CHT] S.~D. Cutkosky, J. Herzog and N.~V. Trung,
Asymptotic behaviour of the Castelnuovo-Mumford regularity,
Composition Math. 118 (1999), 243-261. \par

\noindent [ELS] L. Ein, R. Lazarsfeld and K. Smith, Uniform
bounds and symbolic powers on smooth varieties, Preprint.
\par

\noindent [E] D. Eisenbud, Commutative algebra, with a view
toward algebraic geometry, Springer, 1994. \par

\noindent [EG] D.~Eisenbud and S.~Goto, Linear free resolutions
and minimal multiplicities, J.~Algebra~88 (1984), 89-133.\par

\noindent [GGP] A.~V. Geramita, A. Gimigliano and Y. Pitteloud,
Graded Betti numbers of some embedded rational $n$-folds, Math.
Ann. 301 (1995), no. 2, 363--380. \par

\noindent [GW] S. Goto and K. Watanabe, On graded rings I, J.
Math. Soc. Japan 30 (1978), 179-213. \par

\noindent [HT] L.~T. Hoa and N.~V. Trung, On the
Castelnuovo-Mumford regularity and the arithmetic degree of
monomial ideals, Math. Z. 229 (1998), 519--537. \par

\noindent [K] V. Kodiyalam, Asymptotic behaviour of
Castelnuovo-Mumford regularity, Proc. Amer. Math. Soc. 128
(2000), 407-411. \par

\noindent [NR] D. G. Northcott and D. Rees, Reductions of ideals
in local rings, Proc. Cambridge Philos. Soc. 50 (1954), 145-158.\par

\noindent [SS] K. Smith and I. Swanson, Linear bounds on growth of associated primes, Comm. 
Algebra 25 (1997), 3071--3079.
\par

\noindent [St] B. Sturmfels, Four counterexamples in
combinatorial algebraic geometry, J. Algebra 230 (2000),
282--294. \par

\noindent [Sw] I. Swanson, Powers of ideals. Primary
decompositions, Artin-Rees lemma and regularity, Math. Ann. 307
(1997), 299--313. \par

\noindent [T1] N.~V.~Trung, Reduction exponent and degree bound
for the defining equations of graded rings, Proc. Amer. Math.
Soc.  101 (1987), 229-236. \par

\noindent [T2] N.~V.~Trung,  The Castelnuovo regularity of the
Rees algebra and the associated graded ring, Trans. Amer. Math.
Soc. 350 (1998),  2813--2832. \par

\noindent [T3] N.~V.~Trung, Gr\"obner bases, local cohomology and
reduction number, to appear in Proc. Amer. Math. Soc.\par

\noindent [VV] P.~Valabrega and G.~Valla, Form rings and regular 
sequences, Nagoya Math. J. 72 (1978), 93-101. \par

\noindent [V]  W. Vasconcelos, Cohomological degrees of graded
modules, Six lectures on commutative algebra (Bellaterra, 1996),
345--392, Progr. Math. 166, Birkh\"auser, Basel, 1998. \par

\end{document}